\documentclass[12pt, oneside]{article}   	% use "amsart" instead of "article" for AMSLaTeX format
\usepackage{geometry}                		% See geometry.pdf to learn the layout options. There are lots.
\usepackage{graphicx}				% Use pdf, png, jpg, or eps§ with pdflatex; use eps in DVI mode
\usepackage{amssymb}
\usepackage{amsmath,amsthm,amsfonts,epsfig,amssymb,natbib,eucal,eufrak}

\newtheorem{thm}{Theorem}
\newtheorem{definition}{Definition}%[section]

\title{Brief Article}
\author{The Author}
\date{}							% Activate to display a given date or no date

\title{\bf Calibrating dependence between random elements}
\author{Abram M. Kagan\thanks{Corresponding author}\\
Dept. of Mathematics, University of Maryland, College Park, MD 20742, USA \\
E-mail: amk@math.umd.edu
\and
Gabor J. Sz\'{e}kely\\ 
National Science Foundation, Alexandria, VA 22314 and\\
R\'{e}nyi Institute of Mathematics, Hungarian Acad. of Sciences\\
E-mail:gszekely@nsf.gov}

\begin{document}
\maketitle

\begin{abstract}
\noindent
Attempts to quantify dependence between random elements X and Y via maximal correlation go back to
Gebelein (1941) and R\'{e}nyi (1959).
After summarizing properties (including some new) of the R\'{e}nyi measure of dependence, a calibrated
scale of dependence is introduced. It is based on the ``complexity`` of approximating functions of X by functions of Y.\\
\\
{\it Keywords: Quantification, measures of dependence, projection }
\end{abstract}

\newpage
\section{Introduction}
Let $X,\:Y$  be a pair of random elements taking values in a general probability space \\$\{{\cal X}\times\:{\cal Y},\:{\cal A}\otimes{\cal B},\:P\}$ where $\cal A$ and $\cal B$ are $\sigma$-algebras of subsets of $\cal X$ and $\cal B$, respectively, and $P$ is a probability measure on ${\cal A}\otimes{\cal B}$.\\
The study of dependence of $X$ from $Y$ attracted much attention; see, e. g., Sz\'{e}kely, Rizzo and Bakirov (2007), Sz\'{e}kely and Rizzo (2009), Reshef {\it et al.}(2011), Reimherr and Nicolae (2013). In R\'{e}nyi (1959) natural (``necessary'') requirements to measures of dependence were stated.
Essential features of dependence of $X$ from $Y$ are contained in properties of the conditional expectation
$E\{\varphi(X)|Y\}$ defined for all $\varphi(X)$ with  $E\{|\varphi(X)|\}<\infty.$ If one reduces the domain to $\varphi(X)$ with $E\{|\varphi(X)|^2 \}<\infty$, the conditional expectation becomes the projection of the subspace $L^{2}(X)$ into the subspace of $L^{2}(Y)$ of $L^{2}(X,\:Y)$, the Hilbert space  of functions $\phi(X,\:Y)$ with $E\{|\phi(X,\:Y)|^2\} =\int_{\cal X}\int_{\cal Y} |\phi(x,y)|^2 dP(x,y)<\infty$ and the standard inner product
\[(\phi_1 ,\:\phi_2)=E\{\phi_1 (X,\:Y)\phi_2 (X,\:Y)\}.\]
 According to R\'{e}nyi (1959),in the classical monograph Kolmogorov (1933) introduced the following ratios of variances to measure dependence between $X$ and $Y$:
  \[{\rm var}\:E\{\varphi(X)|Y\}/{\rm var}\:\varphi(X)\]
  Its supremum
\begin{equation}
D(X:Y)=\sup_{||\varphi||=1}{\rm var}\:E\{\varphi(X)|Y\}
\end{equation}
 possesses many properties required from a measure of dependence. These properties are summarized in the following.

\begin{thm}
\begin{itemize}
\item[({\rm i})]$0\leq D(X:Y)\leq 1$.
\item[({\rm ii})]$D(X:Y)=0$ iff $X$ and $Y$ are independent.
\item[({\rm iii})]$\:D(X:Y)=1$ if $X=h(Y)$ (but not iff).
\item[({\rm iv})]$D(X:Y)$ is monotone in the following sense: for any random element\\ $Z,\:D(X:(Y,\:Z))\geq D(X:Y).$
\item[({\rm v})] If $Z$ is independent of the pair $(X,\:Y)$, then $D(X:(Y,\:Z))=D(X:Y)$.
\item[({\rm vi})] For a bivariate Gaussian vector $(X,\:Y)$ with corr$(X,\:Y)=\rho$, $D(X:Y)=|\rho|$.
\item[({\rm vii})] For any bivariate vector $(X,\:Y)$ with corr$(X,\:Y)=\rho$, $D(X:Y)\geq |\rho|$.
\item[({\rm viii})] For an $(m+n)$-dimensional Gaussian random vector $(X,\:Y)$ with variance-covariance matrix
 $$\begin{pmatrix}
 V_{11}& V_{12}\\
 V_{21}& V_{22},
 \end{pmatrix}$$
 assuming the submatrices $V_{11}={\rm var}(X),\:V_{22}={\rm var}(Y)$ non-singular,
  \[D(X:Y)=\sqrt{\lambda(\Sigma)},\]
  where $\lambda(\Sigma)$ is the maximal eigenvalue of the matrix $\Sigma =V_{11}^{-1/2}V_{12}V_{22}^{-1}V_{21}V_{11}^{-1/2}.$
\item[({\rm ix})]
 For an arbitrary $(m+n)$-dimensional vector $(X,\:Y)$ with the above variance-covariance matrix
 \[ D(X:Y)\geq\sqrt{\lambda(\Sigma)}.\]
 \item[({\rm x})]
If a self-decomposable $Z$ is independent of $(X,\:Y)$, then as a function of $\lambda$
 \[D(X:Y+\lambda Z)\: {\rm increases\:on} \:(-\infty.\:0)\: {\rm and\: decreases\:on}\: (0,\infty).\]
 \item[ ({\rm xi})]
 \[D(X:Y)=\{R(X,Y)\}^2\]
 where $R(X,Y)=\sup\rho(\varphi(X),\:\phi(Y)$, $\rho$ is the classical (Pearson) correlation coefficient and the supremum is taken over all $\varphi(X),\:\phi(Y)$ with $0<{\rm var}\varphi(X)<\infty,\\ 0<{\rm var}\phi(Y)<\infty$.
 \item [({\rm xii})] $D(X:Y)=D(Y:X)$.
 \end{itemize}
 \end{thm}
 
 {\it Proofs}.
 Properties (i)-(iv) and (vi)-(ix) are well known, (v) follows from that if $Z$ is independent of $(X,\:Y)$ then for any $\varphi(X)$ with finite expectation
 \[E\{\varphi(X)|Y,\:Z\}=E\{\varphi(X)|Y\}\]
 (see, e. g., Meyer (1966), Theorem 61). Actually, (v) holds if $X$ and $Z$ are conditionally independent given $Y$ which is a weaker condition than independence of $Z$ and $(X,\:Y)$.\\
 To prove (x), recall that a random variable $Z$ is called {\it self-decomposable} if for any $c,\:0<c<1$ there is an independent of $Z$ random variable $Z_c$ such that $Z$ is equidistributed with $cZ + Z_c, \:Z\cong cZ + Z_c$. All random variables having stable distributions are self-decomposable.
 See Lukacs (1970), Ch. 5 for properties of self-decomposable random variables.\\
 Let now $\lambda_1 >\lambda_2 >0,\:\lambda_2 = c\lambda_1$ for some $0<c<1$ and let $Z_c$ be an independent of $X,\:Y,\:Z$ random variable such that $Z\cong cZ + Z_c$. Then

 \begin{eqnarray*}
 & D(X:Y +\lambda_1 Z)=D(X:Y +\lambda_1 (cZ+Z_c))=D(X:Y+\lambda_2 Z +Z_c)\leq D(X:(Y+\lambda_2 Z,Z_c))\\
 & =D(X:Y+\lambda_2 Z)
 \end{eqnarray*}
 due to monotonicity (iv) and (v).\\
 Property (xi) follows from $\rho(\varphi(X),\:\phi(Y))=\rho(E\{\varphi(X)|Y\},\:\phi(Y))$ and the Cauchy-Schwarz inequality
 \[|\rho(\varphi(X),\:\phi(Y))|^2\leq{\rm var}\phi(Y){\rm var}E\{\varphi(X)|Y\}\]
 with the equality sign holding for $\phi(Y)=cE\{\varphi(X)|Y\}$.\\ 
 \\
 Some special forms of dependence are worth studying. One example may be $X,\:Y$ such that
 \begin{equation}
 E\{\varphi(X)|Y\}={\rm const} \Rightarrow P\{\varphi(X)={\rm const}\}=1.
 \end{equation}
 If $p(x|y)$ is the conditional probability density function of $X$ given $Y=y$, (2) is what in statistics is known as {\it completeness} of the family $\{p(x|y),\:y\in{\cal Y}\}$ with $y$ as a family parameter.\\
  Plainly, if $X$ is fully dependent of $Y$, i. e.,$X=h(Y)$, the relation (2) holds.\\
 For independent $X,\:Y,\:  E\{\varphi(X)|Y\}={\rm const}$ for any $\varphi(X)$ with finite expectation, so that (2) moves
 $X,\:Y$ far away from independence.
 \section{The main result}
Here a family $D_{m}(X:Y),\:m=0,\:1,\:2,\ldots$ of measures of dependence is defined that generalizes $D(X:Y)$ in the following direction.\\
A random element $X$  is called $m$-{\it dependent} of $Y$ if for any $\varphi(X)\in L^{2}(X)$, the conditional expectation $E\{\varphi(X)|Y\}=\phi(Y)$ belongs to an $(m+1)$-dimensional subspace of $L^{2}(Y)$.
It is proved that $X$ is $m$-dependent of $Y$ if and only if $D_{m}(X:Y)=0$. For $m=0,\\D_{0}(X:Y)=D(X:Y)$ so that $0$-dependence means independence when for any $\varphi(X),\:E\{\varphi(X)|Y\}$ belongs to the one-dimensional subspace of constants.\\
Denote by $V(\varphi_0,\ldots,\varphi_m)$ the covariance matrix of the $(m+1)$-dimensional vector $(\varphi_{0}(X),\ldots,\varphi_{m}(X))$
of random variables in $L^{2}(X)$ and by $|V|$ the determinant of $V$ (sometimes referred to as the generalized variance of the vector with covariance matrix $V$.\\
For $\varphi(X)\in L^2 (X)$ set $\phi(Y)=E\{\varphi(X)|Y\}$.
\begin{definition} The index of $m$-dependence of $X$ from $Y$ is defined as
\begin{equation}
D_{m}(X:Y)=\sup|V(\phi_0 (Y),\phi_1 (Y),\ldots,\phi_m (Y))|,
\end{equation}
where the supremum is taken over all $\varphi_0 (X),\varphi_1 (X),\ldots,\varphi_m (X)$ with\\
$V(\varphi_0,\varphi_1,\ldots,\varphi_m)={\rm Id}_{m+1}$, the identity matrix of order $m+1$.
\end{definition}
The matrix $V(\phi_0 (Y),\phi_1 (Y),\ldots,\phi_m (Y))$ is non-negative definite and, as well known,
\[V(\phi_0 (Y),\phi_1 (Y),\ldots,\phi_m (Y))\leq V(\varphi_0,\varphi_1,\ldots,\varphi_m),\]
so that $0\leq D_{m}(X:Y)\leq 1.$\\
Plainly, $D_{m}(X:Y)$ is monotone in the same sense as $D(X:Y)$ (see (iv)).\\
Since $D_{0}(X:Y)=\sup_{\varphi:{\rm var}\varphi(X)=1}{\rm var}E\{\varphi(X)|Y\}$, $D_{0}(X:Y)=0$ if and only if
for any $\varphi\in L^2 (X)$, $E(\varphi|Y)$ belongs to the one-dimensional subspace of constants in $L^2 (Y)$. A similar property of $D_{m}(X:Y)$ is proved in the next theorem.
\begin{thm}
$D_{m}(X:Y)=0$ if and only if for any $\varphi\in L^2 (X)$, $E(\varphi|Y)$ belongs to a subspace of $L^2 (Y)$ of dimension $\leq m+1$.\\
\end{thm}
{\it Proof}. First notice that the image of $\varphi \equiv 1\in L^2 (X)$ is
$\phi\equiv 1\in L^2(Y)$. For any $\varphi\in L^2(X)$ with $E(\varphi(X))=0$ its image
$E(\varphi|Y)=\phi(Y)$ also has $E(\phi(Y))=0$.
We shall show that $D_{m}(X:Y)=0$ if and only if for any $\varphi\in L^2 (X)$,  with $E(\varphi)=0$ its image $E(\varphi|Y)$ belongs to a subspace of $L^2 (Y)$ of dimension $\leq m.$\\
Now take elements $\varphi_1,\ldots,\varphi_{m+1}$ in $L^2 (X)$ with zero expectations and the identity covariance matrix (that  is why $\varphi \equiv 1$ is treated separately).\\ Due to $D_m(X:Y)=0$, the covariance matrix $V(\phi_1,\ldots,\phi_{m+1})$ of their images $\phi_1 (Y)=E(\varphi_1 |Y),\ldots,\phi_{m+1} (Y)=E(\varphi_{m+1} |Y)$ is degenerate and thus
 \[c_1 \phi_0 +\ldots+c_{m+1} \phi_m =0\]
 with probability 1 for some constants $c_1,\ldots,c_{m+1}$.\\
 Assuming without loss in generality $\phi_1,\ldots,\phi_{m}$ linearly independent and thus $c_{m+1}\neq 0$, $D_{m}(X:Y)=0$ implies
 that for the chosen $\varphi_1,\ldots,\varphi_{m+1}$ their images belong to the $m$-dimensional subspace span$\{\phi_1,\ldots,\phi_{m}\}$. If $\varphi\in{\rm span}\{\varphi_1,\ldots,\varphi_{m+1}\}$, then plainly
 $E(\varphi|Y)=\phi\in{\rm span}\{\phi_1,\ldots,\phi_{m}\}$.\\
 Let now $\varphi\in L^2 (X)$ belong to the orthocomplement of span$\{\varphi_1,\ldots,\varphi_{m+1}\},
 \varphi\in({\rm span}\{\varphi_1,\ldots,\varphi_{m+1}\})^{\perp}$. Assuming var$\varphi(X)=1$,
 consider the vector $(\varphi_1,\ldots,\varphi_{m},\varphi)$. The covariance matrix $V(\phi_1,\ldots,\phi_{m},\phi)$ where $\phi=E(\varphi|Y)$, is degenerate and for some constants $c'_1,\ldots,c'_{m+1}$
 \[c'_1\phi_1+\ldots+c'_{m}\phi_{m}+c'_{m+1}\phi_m\]
 so that $\phi\in{\rm span}\{\phi_1,\ldots,\phi_{m}\}$.\\
 Since any $\varphi\in L^2 (X)$ can be represented as
 \[\varphi = \varphi' + \varphi'',\:\varphi' \in {\rm span}\{\phi_1,\ldots,\phi_{m+1}\},\:\varphi'' \in ({\rm span}\{\phi_1,\ldots,\phi_{m+1}\})^\perp,\]
 the above proves the sufficiency part of Theorem 2.\\
 To prove necessity, take any set $\varphi _1,\ldots,\varphi_{m+1}$ of elements of $L^2 (Y)$ with zero expectations and identity covariance matrix. If their images
 \[  \phi_1 =E(\varphi_1 |Y),\ldots, \phi_{m+1}=E(\varphi_{m+1}|Y)\]
 belong to an $m$-dimensional subspace of $L^2 (Y)$ and thus are linearly dependent,
  their covariance matrix is degenerate and $|V(\phi_1,\ldots,\phi_{m+1})|=0.$\\
  The images of constants form an one-dimensional subspace of constants and for any $\varphi(X),\:\varphi = E(\varphi)+[\varphi -E(\varphi)]$ proving the necessity part of
  Theorem 2. \\
  A trivial corollary is $D_{m}(X:Y)=0\Rightarrow D_{m+1}(X:Y)=0$.\\
 Starting with the class ${\cal C}_0$ of independent $X,\:Y$, the indices $D_m$ lead to a scale of classes ${\cal C}_m$,
 \[{\cal C}_0 \supset{\cal C}_1\supset{\cal C}_2\supset\ldots,\]
 where $(X,\:Y)\in{\cal C}_m \Leftrightarrow D_{m}(X:Y)=0$. \\
 The classes ${\cal C}_m$ are non-empty. If the conditional probability density function $p(x|y)$ of $X$ given $Y$ has a form of
 \[p(x|y)=p_0 (x) +p_1 (x)q_1 (y)+\ldots+p_m (x)q_m (y),\]
 with $q_1 \in L^2 (Y),\ldots, q_m \in L^2 (Y)$, then
 \[E(\varphi (X)|Y)=\int \varphi (x) p(x|Y)dx\]
 belongs to the $(m+1)$-dimensional subspace span $\{1,q_{1}(Y),\ldots,q_{m}(Y)\}$ of $L^2 (Y)$,
 so that $(X,\:Y)\in {\cal C}_m$.\\
 \\
 Of some interest may be an approach that starts not at independent $X,\:Y$ but at the other end when $X$ and $Y$ are arbitrarily dependent. Let us call $X$ $m$-codependent of $Y$ if for any $\varphi\in L^2 (X),\:E(\varphi |Y)$ belongs to a subspace of $L^2 (Y)$ of 
 codimension $m$. As $m$ increases, $X$ and $Y$ become less and less dependent, in a sense. How to quantify the codependence?
 \newpage
 \section*{References}
 Gebelein, H. (1941). Das  statistische  Problem  der  Korrelation  als  Variations-  und  Eigenwert-problem  und  sein  Zusammenhang  mit der Ausgleichungsrechnung.
{\it Z. Angew. Math. Mech.}, {\bf 21}, 364--379.\\
\\
 Kolmogoroff, A. N. (1933). {\it Grundbegriffe der Wahrscheinlichkeitsrechnung}, Springer, Berlin.\\
 \\
 Lukacs, E. (1970).{\it Characteristic Functions}, Griffin, London.\\
 \\
 Meyer, P. A. (1966) {\it Probability and Potentials}, Blaisdel.\\
 \\
 Reimherr, R., and Nicolae, D. L. (2013). On quantifying dependence: a framework for developing interpretable measures. {\it Statistical Science}, {\bf 28}, 116-130.\\
 \\
 R\'{e}nyi, A. (1959). On measures of dependence. {\it Acta Math. Acad. Sci. }, {\bf 10}, 441-451.\\
 \\
Reshef, D. N., Reshef, Y. A., Finucane, H. K., Grossman, S. R., McVean, G., Turnbaugh, P. J., Lander, E. S., Mitzenmacher, M., Sabeti, P. C. (2011). Detecting novel associations in large data sets. {\it Science}, {\bf 334}, 1518-1524.\\
\\
Sz\'{e}kely, G. J., Rizzo, M. L., and Bakirov, N. K. (2007).Measuring and testing dependence by correlation of distances. {\it Ann. Statist.}, {\bf 35},2769-2794. \\
\\
Sz\'{e}kely, G. J., and Rizzo, M. L. (2009). Brownian distance covariance. {\it Ann. Appl. Stat.},{\bf 3}, 1236-1265. 
 \end{document}